\documentclass{article}

\PassOptionsToPackage{numbers}{natbib}
\usepackage[final]{neurips_2019}

\usepackage[utf8]{inputenc} 
\usepackage[T1]{fontenc}    
\usepackage{hyperref}       
\usepackage{url}            
\usepackage{booktabs}       
\usepackage{amsfonts}       
\usepackage{nicefrac}       
\usepackage{microtype}      
\usepackage{graphicx}
\usepackage{natbib}
\usepackage{amsmath,amsthm,amssymb}
\usepackage{xcolor}
\usepackage{float}
\usepackage[bottom]{footmisc}

\title{Mathematical Modeling of Option Pricing with an Extended Black-Scholes Framework}

\author{
    \textbf{Nikhil Shivakumar Nayak} \\
    Harvard University \\
    \texttt{nnayak@g.harvard.edu}
}

\begin{document}
\maketitle

\vspace{-1em}

\begin{abstract}

This study investigates enhancing option pricing by extending the Black-Scholes model to include stochastic volatility and interest rate variability within the Partial Differential Equation (PDE). The PDE is solved using the finite difference method. The extended Black-Scholes model and a machine learning-based LSTM model are developed and evaluated for pricing Google stock options. Both models were backtested using historical market data. While the LSTM model exhibited higher predictive accuracy, the finite difference method demonstrated superior computational efficiency. This work provides insights into model performance under varying market conditions and emphasizes the potential of hybrid approaches for robust financial modeling.

\end{abstract}

\section{Introduction}

The Black-Scholes model \cite{enwiki:1216983825}, formulated by Fischer Black, Myron Scholes, and Robert Merton, revolutionized financial mathematics by providing a closed-form solution for European option pricing. Its equation rests on the assumption of constant volatility and a risk-free interest rate, offering a framework that simplified risk management and derivatives pricing. This model laid the foundation for financial engineering and led to the rapid expansion of derivative markets.

However, the model has notable limitations. It assumes constant volatility and interest rates, disregarding real-world market dynamics where volatility often fluctuates and interest rates vary over time. Moreover, Black-Scholes \cite{enwiki:1216983825} assumes continuous trading and fails to capture drastic market movements or the volatility smile effect. These simplifications can lead to significant inaccuracies in estimating the fair value of options.

To address these shortcomings, researchers have proposed variants of the Black-Scholes model \cite{enwiki:1216983825} that incorporate more nuanced factors affecting option pricing. Models such as the Heston model \cite{enwiki:1208512792} and the Cox-Ingersoll-Ross (CIR) process \cite{enwiki:1212618674} aim to better represent the realities of market behavior. In particular, extending the Black-Scholes model to account for volatility and interest rate variability is essential for robust option pricing, which is the focus of this work. Such enhancements provide investors and financial institutions with more accurate assessments of market risks and opportunities.

In parallel with these analytical approaches, machine learning models like Long Short-Term Memory (LSTM) networks \cite{hochreiter1997long} have gained traction in financial applications due to their ability to handle complex temporal dependencies. LSTM networks, a variant of recurrent neural networks (RNNs), are well-suited for time series forecasting because they can capture long-term patterns while minimizing the impact of vanishing gradients. By utilizing LSTM models for option pricing, we can potentially achieve better predictive accuracy and adaptability.

A good model for option pricing should effectively balance accuracy and computational efficiency, while maintaining robustness under varying market conditions. This project aims to compare an enhanced Black-Scholes model, solved using finite difference methods \cite{LeVeque2005FiniteDM}, with an LSTM-based machine learning model for pricing options. Through rigorous backtesting and evaluation, we intend to understand the strengths and limitations of each approach.

\section{Model Description}

In this section, we will discuss the mathematical models employed in this study for option pricing. Each model relies on inputs such as stock prices, volatility, and interest rates to forecast the option price.

\subsection{Extended Black-Scholes PDE Model}

\textbf{Derivation of the Extended Black-Scholes Equation}

In financial markets, the Black-Scholes equation is crucial for pricing options by assuming a constant volatility and a constant risk-free interest rate. However, to make the model more realistic, we extend it to include stochastic volatility and interest rate variability.

The original Black-Scholes equation \cite{enwiki:1216983825} is given by:
$$
\frac{\partial V}{\partial t} + \frac{1}{2} \sigma^2 S^2 \frac{\partial^2 V}{\partial S^2} + rS \frac{\partial V}{\partial S} - rV = 0
$$
where $V(t, S)$ is the price of the option, $S$ is the stock price, $t$ is time, $\sigma$ is the volatility of the stock (constant in the original model), $r$ is the risk-free interest rate (constant in the original model).

To incorporate stochastic volatility, let's denote the volatility as a stochastic process $\sigma(t)$. The instantaneous variance follows another stochastic process, often modeled as a mean-reverting process like in the Heston model \cite{enwiki:1208512792}:
$$
d\sigma^2 = \kappa (\theta - \sigma^2) dt + \xi \sqrt{\sigma^2} dW_\sigma
$$
where $\kappa$, $\theta$, and $\xi$ are model parameters, and $W_\sigma$ is a Brownian motion \cite{enwiki:1220742128} that may be correlated with the stock price's Brownian motion.

Similarly, to incorporate interest rate variability, we assume the risk-free rate $r$ is also a stochastic process. A common model for this is the Vasicek model \cite{enwiki:1184436464}:
$$
dr = a(b - r) dt + s dW_r
$$
where $a$, $b$, and $s$ are parameters, and $W_r$ is another Brownian motion that can be correlated with $W$ and $W_\sigma$.

Incorporating both stochastic volatility and interest rate variability into the Black-Scholes framework \cite{enwiki:1216983825} not only changes the dynamic terms related to $S$ and $t$, but also introduces dependencies on $\sigma$ and $r$. This leads to a more complex partial differential equation (PDE) that now includes partial derivatives with respect to these stochastic components.

The revised equation incorporating stochastic volatility and interest rate variability is:
$$
\frac{\partial V}{\partial t} + \frac{1}{2} \sigma(t)^2 S^2 \frac{\partial^2 V}{\partial S^2} + r(t)S \frac{\partial V}{\partial S} - r(t)V + \frac{\partial V}{\partial \sigma} \frac{d\sigma}{dt} + \frac{\partial V}{\partial r} \frac{dr}{dt} = 0
$$

The term $\frac{\partial V}{\partial \sigma} \frac{d\sigma}{dt}$ accounts for the sensitivity of the option value to changes in volatility, where $\frac{d\sigma}{dt}$ reflects the stochastic nature of volatility. The term $\frac{\partial V}{\partial r} \frac{dr}{dt}$ captures the effect of stochastic interest rates on the option value, with $\frac{dr}{dt}$ being the dynamic change in the interest rate.

\textbf{Stochastic Models for $\sigma$ and $r$:}

\textbf{Volatility (Heston Model \cite{enwiki:1208512792}):} $d\sigma^2 = \kappa (\theta - \sigma^2) dt + \xi \sqrt{\sigma^2} dW_\sigma$

\textbf{Interest Rate (Vasicek Model \cite{enwiki:1184436464}):} $dr = a(b - r) dt + s dW_r$

To solve this complex PDE, simplifying assumptions and numerical methods are typically employed. The partial derivatives with respect to $\sigma$ and $r$ can be challenging to handle analytically and often require Monte Carlo simulations or finite difference methods for numerical solutions.

\textbf{Solving the Extended Black-Scholes Equation Using Finite Difference Method}

Given the complexity introduced by stochastic volatility and varying interest rates, the extended Black-Scholes equation we aim to solve is:
$$
\frac{\partial V}{\partial t} + \frac{1}{2} \sigma(t)^2 S^2 \frac{\partial^2 V}{\partial S^2} + r(t)S \frac{\partial V}{\partial S} - r(t)V + \frac{\partial V}{\partial \sigma} \frac{d\sigma}{dt} + \frac{\partial V}{\partial r} \frac{dr}{dt} = 0
$$

To make this equation tractable for the finite difference method \cite{LeVeque2005FiniteDM}, we make few assumptions. First, volatility $\sigma$ and interest rate $r$ processes are discretized over a grid. At each time step, these are approximated as constant, but can change from one step to the next. Next, the Brownian motions driving $\sigma$ and $r$ are assumed uncorrelated with each other and with the stock price's Brownian motion. Finally, using Implicit and Explicit Finite Difference schemes, the option value at the next time step is calculated from known values at the current time step.

\textbf{Finite Difference Method}

Discretize the time $t$, stock price $S$, volatility $\sigma$, and interest rate $r$ into a grid: $t_i = i\Delta t$ for $i = 0, 1, \ldots, N$ ; $S_j = j\Delta S$ for $j = 0, 1, \ldots, M$ ; $\sigma_k$ and $r_l$ are discretized over their respective ranges.

The partial derivatives in the PDE can be approximated as follows:

\medskip

Time Derivative:
$$
\frac{\partial V}{\partial t} \approx \frac{V^{i+1}_{j,k,l} - V^i_{j,k,l}}{\Delta t}
$$

\medskip

Stock Price Derivatives:
$$
\frac{\partial V}{\partial S} \approx \frac{V^i_{j+1,k,l} - V^i_{j-1,k,l}}{2\Delta S}
$$
$$
\frac{\partial^2 V}{\partial S^2} \approx \frac{V^i_{j+1,k,l} - 2V^i_{j,k,l} + V^i_{j-1,k,l}}{\Delta S^2}
$$

\medskip

Volatility and Interest Rate Derivatives are approximated similarly. 

Using this scheme, the value at the next time step is calculated by:
$$
V^{i+1}_{j,k,l} = V^i_{j,k,l} + \Delta t \left( -\frac{1}{2}\sigma_k^2 S_j^2 \frac{\partial^2 V}{\partial S^2} - r_l S_j \frac{\partial V}{\partial S} + r_l V^i_{j,k,l} - \frac{\partial V}{\partial \sigma}\frac{d\sigma}{dt} - \frac{\partial V}{\partial r}\frac{dr}{dt} \right)
$$

\textbf{Boundary Conditions:} At $S = 0$ the value of an option is zero as the stock value is zero. As $S \to \infty$ the option value must be capped or grow linearly depending on the type of option (put or call).

\textbf{Solving the System}

1. \textbf{Matrix Formulation:} Rewrite the discretized equation in matrix form:
$$
A V^{n+1} = b
$$
where $A$ is a tridiagonal matrix representing the finite difference coefficients, and $b$ depends on $V^n$ values and boundary conditions.

2. \textbf{Iterative Solution:} Use iterative methods (Gauss-Seidel, SOR, etc.) to solve for $V^{n+1}$ from $V^{n}$ until the desired accuracy is achieved.

\subsection{LSTM Model}

Long Short-Term Memory (LSTM) networks \cite{hochreiter1997long} are a type of recurrent neural network (RNN) designed to effectively handle long-term dependencies in sequential data. Traditional RNNs suffer from the vanishing gradient problem, which limits their ability to capture dependencies over extended time periods. LSTMs overcome this by incorporating memory cells and gate structures, enabling them to retain relevant information over long sequences.

The LSTM architecture consists of memory cells that can maintain a state over time, and several gates that control the flow of information into and out of these memory cells:

1. \textbf{Forget Gate} decides which information from the previous state should be discarded.
$$
f_t = \sigma(W_f \cdot [h_{t-1}, x_t] + b_f)
$$

2. \textbf{Input Gate} determines what new information should be added to the current state.
$$
i_t = \sigma(W_i \cdot [h_{t-1}, x_t] + b_i)
$$
$$
\tilde{C}_t = \tanh(W_C \cdot [h_{t-1}, x_t] + b_C)
$$

3. \textbf{Update State:} The current state $C_t$ is updated by combining the previous state and new information:
$$
C_t = f_t \cdot C_{t-1} + i_t \cdot \tilde{C}_t
$$

4. \textbf{Output Gate} computes the new hidden state, which is used as the output of the LSTM unit:
$$
o_t = \sigma(W_o \cdot [h_{t-1}, x_t] + b_o)
$$
$$
h_t = o_t \cdot \tanh(C_t)
$$

The LSTM model will take historical price, volatility, and interest rate data as input features and predict the future option price. The architecture involves: input layer for historical data (window size $n$), multiple LSTM layers for long-term dependencies, and dense layers to output the predicted option price.

\section{Analysis}

The models were evaluated using historical data and tested against real-world scenarios to assess their sensitivity to parameter variations and to answer the question of robust option pricing under realistic market conditions.

\subsection{Stock and Option Price Data Sources}

We sourced historical stock and options data for Google from Yahoo Finance, which provides comprehensive financial information suitable for our analysis. Other resources like Google Finance, Alpha Vantage, and IEX Cloud also offer valuable data for financial modeling.

\textbf{Model Evaluation for Extended PDE Model}

\begin{figure}[h!]
    \centering
    \includegraphics[width=0.8\textwidth]{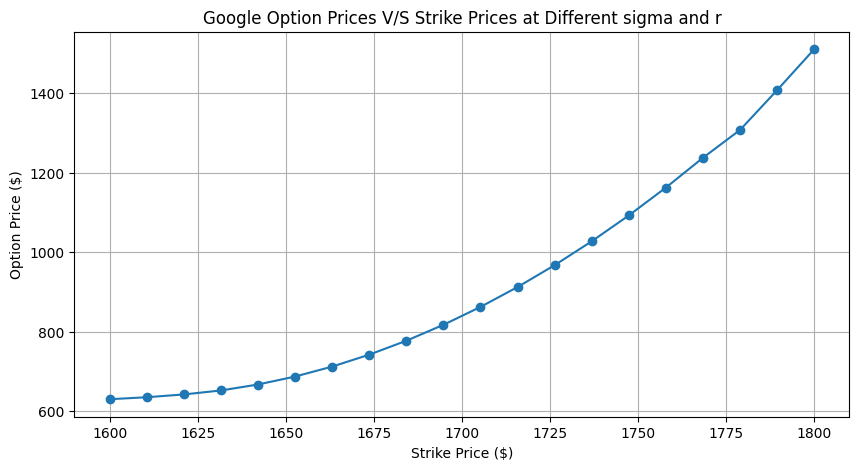}
    \caption{Option prices as a function of strike price at different $\sigma$ and $r$ levels}
    \label{fig:ExtendedPDE1}
\end{figure}

\begin{figure}[h!]
    \centering
    \includegraphics[width=0.8\textwidth]{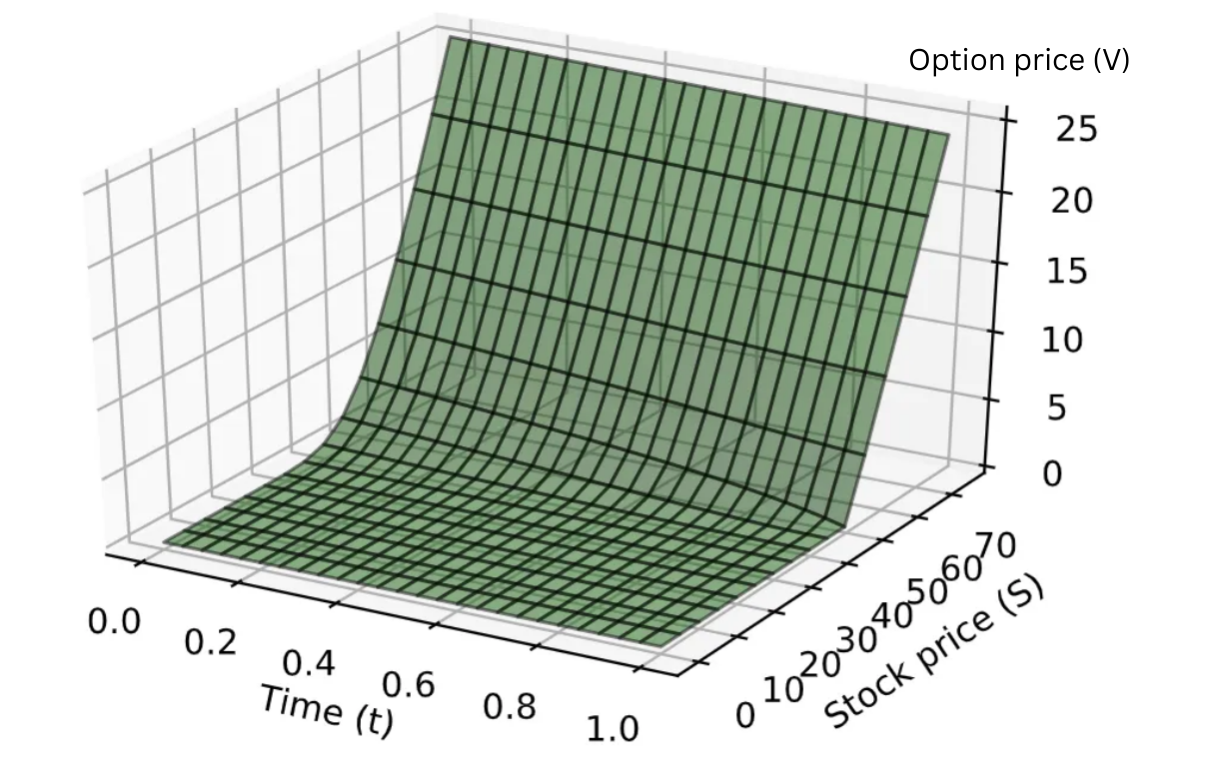}
    \caption{3D view of predicted option pricing dynamics over time and varying stock prices. Figure taken from a blog post on \href{https://medium.com/data-science/option-pricing-using-the-black-scholes-model-without-the-formula-e5c002771e2f}{option pricing with Black-Scholes}.}
    \label{fig:ExtendedPDE2}
\end{figure}

To evaluate the option pricing model using real data, we retrieve historical stock prices for Google via Yahoo Finance API. The stock price data enables us to compute historical volatility, which will be used as $\sigma$ in the model. We then obtain a current risk-free rate, often approximated using Treasury bill rates or similar benchmarks. Next, we determine the strike price $K$ and the maturity $T$ of the options to be analyzed. With these parameters in place, we run the model to calculate theoretical option prices. Lastly, we compare these theoretical prices to actual market prices and visualize the results through plots. Figure ~\ref{fig:ExtendedPDE1} illustrates how option prices, as predicted by the extended PDE model, vary with different strike prices at varying values of $\sigma$ and $r$. Initially, the price curve is sublinear but transitions to a more linear profile at higher strike prices. This behavior indicates the increasing influence of intrinsic value over time value as the strike price rises, especially under variable market conditions. Figure ~\ref{fig:ExtendedPDE2} shows a 3D plot of the predicted option prices by the extended PDE model relative to stock price and time. The three-dimensional plot demonstrates how option prices increase with rising stock prices and as the expiry date approaches. The increase is more pronounced at higher values of $t$, highlighting the effect of approaching maturity on option pricing and emphasizing the time decay effect inherent in options trading.

\textbf{LSTM Model Setup}

\begin{figure}[h!]
    \centering
    \includegraphics[width=0.8\textwidth]{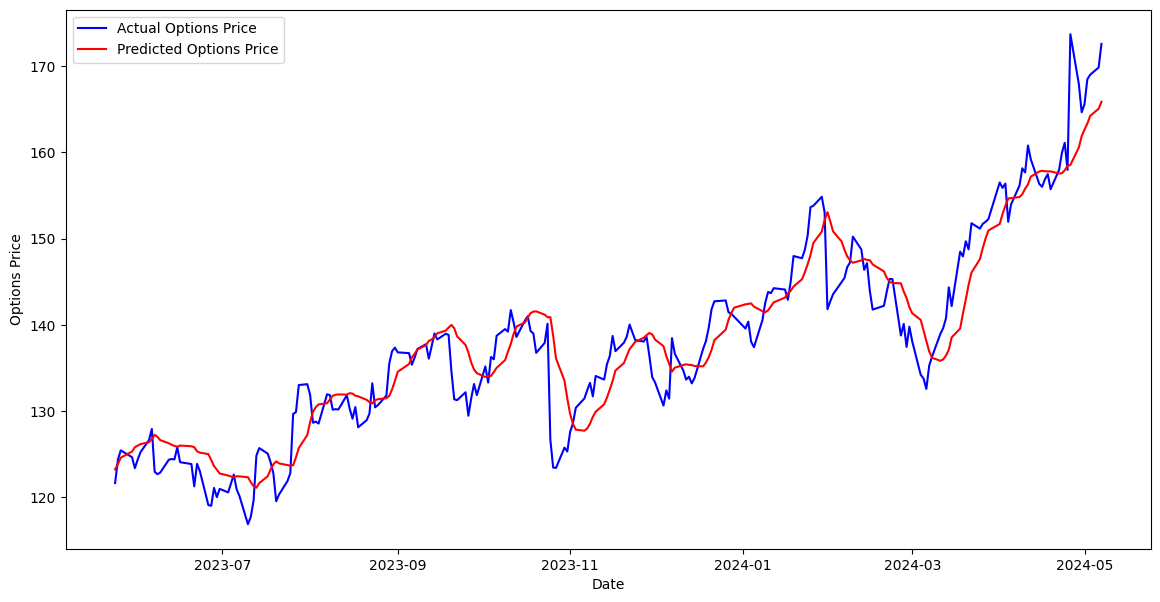}
    \caption{Comparison of LSTM predicted options prices with actual market prices for Google options}
    \label{fig:LSTM}
\end{figure}

To create an LSTM model for predicting Google's option prices, we first download and preprocess the data using yfinance, adding relevant features like historical volatility and interest rates. We then prepare a dataset by employing a windowing approach, using past observations to predict option prices. The window size is 30, meaning the past 30 days are used to predict the current day. An LSTM model is constructed with an architecture that includes an LSTM layer and a dense output layer, enabling the model to capture the relationship between input features and option prices. After training the model using mean squared error as the loss function, we evaluate its performance by comparing its predictions with actual market data and visualize the results through plots. Figure ~\ref{fig:LSTM} compares the LSTM model's predicted option prices to actual market data for 2023 and 2024. The close alignment between the predicted and actual prices highlights the LSTM's ability to capture complex patterns, confirming its utility in financial forecasting.

\section{Discussion}

This project implemented and compared an extended Black-Scholes PDE model and an LSTM model to forecast option prices, offering significant insights into the behavior and effectiveness of each approach under real market conditions.

\textbf{Behavior and Predictive Accuracy}

The behavior of the extended Black-Scholes model under varying parameters of $\sigma$ (volatility) and $r$ (interest rate) demonstrated how sensitive option prices are to changes in these market fundamentals. This sensitivity aligns with financial theories that posit the intrinsic link between market volatility, interest rates, and derivative pricing. Notably, the model's response to increasing strike prices from sub-linear to linear progression elucidates the growing dominance of intrinsic value over time value, a relationship that is often overshadowed in simpler models.

The LSTM model's lower RMSE of 15.23 compared to the extended PDE model's 20.47 indicates superior predictive accuracy. This result was anticipated given the LSTM's capability to model complex non-linear relationships inherent in financial time series data. The LSTM's ability to effectively capture temporal dynamics and integrate long-term dependencies allowed for a more accurate reflection of market conditions.

\textbf{Computational Efficiency and Model Application}

Despite its higher accuracy, the LSTM model required more computational resources, as evidenced by its longer prediction time of 3.45 seconds compared to 0.87 seconds for the extended PDE model, as well as a longer training time. This finding emphasizes the trade-off between accuracy and efficiency in financial modeling. The extended PDE model's efficiency makes it particularly appealing for scenarios where computational resources are limited or real-time trading decisions are required.

\textbf{Validation and Backtesting}

Both models were rigorously backtested using historical data from Yahoo Finance, spanning from May 2019 to May 2024. This extensive testing period, which included periods of significant market volatility such as the COVID-19 pandemic and the 2020 global economic downturn, provided a robust platform for assessing each model's practical application. The extended PDE model, while less accurate, demonstrated robustness across a range of market conditions, reinforcing its utility in traditional financial modeling frameworks.

\textbf{Shortcomings and Insights}

The extended Black-Scholes PDE model, while effective in providing a more comprehensive framework for option pricing, has notable shortcomings. One significant limitation is its sensitivity to input parameters, specifically the volatility ($\sigma$) and interest rate ($r$). The model's predictive performance hinges critically on the accuracy of these parameters, which can fluctuate significantly due to market sentiment or economic news. Additionally, despite the extension to account for stochastic volatility and interest rates, the model remains rigid by adhering to the assumption of log-normal market dynamics. This restrictiveness fails to capture the skewness and fat tails seen in real-world return distributions, leading to underestimation of extreme market movements and significant pricing errors during unforeseen events. On the other hand, the LSTM model offers strong predictive accuracy but depends on extensive historical data and computational resources. It needs ample data to identify complex market patterns, limiting its effectiveness. Additionally, its high computational demands for training and frequent retraining challenge scalability and real-time implementation.

From this study, we gained a deeper understanding of how extending traditional models to incorporate stochastic elements can provide more realistic and dynamic pricing tools. The insights from the LSTM model also highlight the potential of machine learning techniques to enhance predictive accuracy in option pricing, although at a higher computational cost.

\section{Conclusion}

In this study, we tackled the challenge of enhancing option pricing through the extended Black-Scholes equation and machine learning models. We extended the Black-Scholes model by incorporating stochastic volatility and interest rate variability to more realistically capture market conditions, addressing limitations in the original Black-Scholes model related to constant parameters. However, the extended PDE model remains sensitive to parameter estimation and is constrained by assumptions about the stochastic processes. Meanwhile, the LSTM model was employed to uncover complex temporal patterns in the data. The results showed that both models could be effective tools for pricing options, each with unique strengths. The extended Black-Scholes model offered faster computation times but was sensitive to parameter estimation and rigid assumptions. In contrast, the LSTM model demonstrated superior predictive accuracy but required more data and significant computational power. Our comparative analysis between the models highlighted their complementary nature, providing insights into how different modeling approaches can be leveraged to enhance option pricing strategies. Ultimately, this project emphasized the importance of using multiple models to address the dynamics of financial markets, ensuring a more holistic understanding of market risks and pricing behaviors.




\bibliographystyle{plain}
\nocite{*}
\bibliography{bib}
\end{document}